\documentclass[12pt]{article}
\usepackage{amssymb}

\begin{document}

\title{The velocity translation in the game of ''Life''}
\author{Ed Troyan \\
\textit{Academy for Management of Innovations }\\
\textit{16a Novobasmannaya Street, Moscow 107078, }\\
\textit{Russian Federation}}
\maketitle

\begin{abstract}
The velocity translation law for the patterns in the game of 'Life'  is
discovered. It is found for an arbitrary angle between the velocities in the
moving reference frames. The formula differs from its physical prototype of
special relativity but admits gentle reduction to the Galilean limit and
provides absence of extra-light speeds.
\end{abstract}

\sloppy\tabcolsep=0.15cm

\section{Introduction}

The game of \textsc{Life} is the most famous cellular automata introduced by
J. Conway \cite{r1,r2}. In 1970 it appeared as a mere recreation puzzle, now
it has found applications in various scopes of science, from pure game
theory to biology, chemistry, economics etc. The natural and simple rules
make this game popular for the wider community also. In the initial
formulation of Conway it acts on a two-dimensional square lattice where the
cells exist in two states -- dead (0) or alive (1). The evolution of each
cell depends on the status of its eight neighbors.

A dead cell comes to live, if it has exactly \textit{three }living neighbors.

A cell remains living if it has \textit{two} or \textit{three} living
neighbors; otherwise (without support) it will die.

The name of this cellular automata reflects a colony of primitive one-cell
microorganisms: both rare and dense population has no chance to survive,
while the intermediate densities may admit intensive growth and development.
This seeming simplicity is accompanied by extremely complicated dynamics,
investigated in the recent years.

One of the most striking feature of this game of \textsc{Life} is the
existence of cell societies endowed with definite regularity in their
evolution. The so-called \textit{spaceships} form the most important class
of such objects. \textit{Spaceship} is the pattern which reproduces itself
after $P$ moves ('period') and appears \textit{replaced} in several squares,
vertical or horizontal, or in both directions. If $n$ is the magnitude of
this replacement the quantity $v=n/P$ implies the velocity of the spaceship.
'Glider (below) is the smallest spaceship 
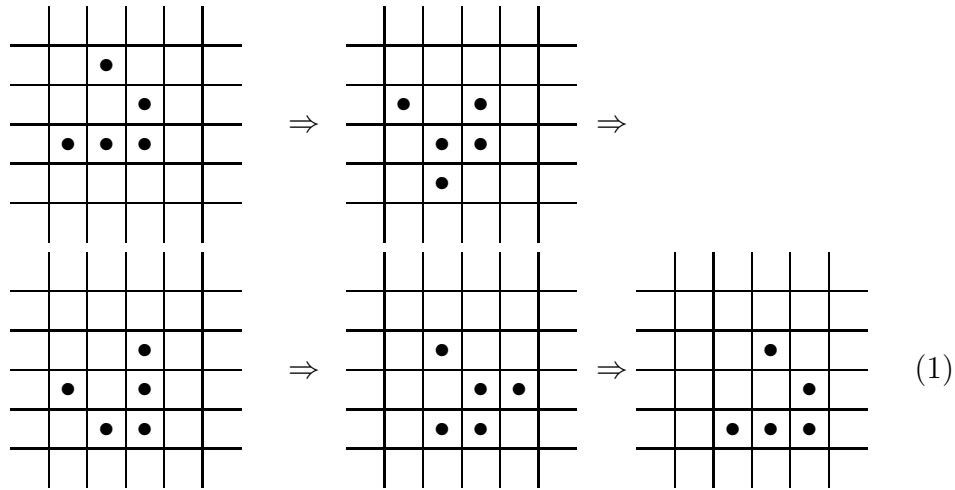
\begin{eqnarray}
\begin{tabular}{c|c|c|c|c|c}
\phantom{$\bullet$} &  &  &  &  &  \\ \hline
&  & $\bullet $ &  &  &  \\ \hline
&  &  & $\bullet $ &  &  \\ \hline
& $\bullet $ & $\bullet $ & $\bullet $ &  &  \\ \hline
&  &  &  & \phantom{$\bullet$} &  \\ \hline
&  &  &  &  & \phantom{$\bullet$}%
\end{tabular}%
\ \  &\Rightarrow &%
\begin{tabular}{cccccc}
\multicolumn{1}{c|}{\phantom{$\bullet$}} & \multicolumn{1}{c|}{} & 
\multicolumn{1}{c|}{} & \multicolumn{1}{c|}{} & \multicolumn{1}{c|}{} &  \\ 
\hline
\multicolumn{1}{c|}{} & \multicolumn{1}{c|}{} & \multicolumn{1}{c|}{} & 
\multicolumn{1}{c|}{} & \multicolumn{1}{c|}{} &  \\ \hline
\multicolumn{1}{c|}{} & \multicolumn{1}{c|}{$\bullet $} & 
\multicolumn{1}{c|}{} & \multicolumn{1}{c|}{$\bullet $} & 
\multicolumn{1}{c|}{} &  \\ \hline
\multicolumn{1}{c|}{} & \multicolumn{1}{c|}{} & \multicolumn{1}{c|}{$\bullet 
$} & \multicolumn{1}{c|}{$\bullet $} & \multicolumn{1}{c|}{} &  \\ \hline
\multicolumn{1}{c|}{} & \multicolumn{1}{c|}{} & \multicolumn{1}{c|}{$\bullet 
$} & \multicolumn{1}{c|}{} & \multicolumn{1}{c|}{\phantom{$\bullet$}} &  \\ 
\hline
\multicolumn{1}{c|}{} & \multicolumn{1}{c|}{} & \multicolumn{1}{c|}{} & 
\multicolumn{1}{c|}{} & \multicolumn{1}{c|}{} & \phantom{$\bullet$}%
\end{tabular}%
\ \Rightarrow   \nonumber  \label{g} \\
\begin{tabular}{cccccc}
\multicolumn{1}{c|}{\phantom{$\bullet$}} & \multicolumn{1}{c|}{} & 
\multicolumn{1}{c|}{} & \multicolumn{1}{c|}{} & \multicolumn{1}{c|}{} &  \\ 
\hline
\multicolumn{1}{c|}{} & \multicolumn{1}{c|}{} & \multicolumn{1}{c|}{} & 
\multicolumn{1}{c|}{} & \multicolumn{1}{c|}{} &  \\ \hline
\multicolumn{1}{c|}{} & \multicolumn{1}{c|}{} & \multicolumn{1}{c|}{} & 
\multicolumn{1}{c|}{$\bullet $} & \multicolumn{1}{c|}{} &  \\ \hline
\multicolumn{1}{c|}{} & \multicolumn{1}{c|}{$\bullet $} & 
\multicolumn{1}{c|}{} & \multicolumn{1}{c|}{$\bullet $} & 
\multicolumn{1}{c|}{} &  \\ \hline
\multicolumn{1}{c|}{} & \multicolumn{1}{c|}{} & \multicolumn{1}{c|}{$\bullet 
$} & \multicolumn{1}{c|}{$\bullet $} & \multicolumn{1}{c|}{%
\phantom{$\bullet$}} &  \\ \hline
\multicolumn{1}{c|}{} & \multicolumn{1}{c|}{} & \multicolumn{1}{c|}{} & 
\multicolumn{1}{c|}{} & \multicolumn{1}{c|}{} & \phantom{$\bullet$}%
\end{tabular}%
\ \  &\Rightarrow &%
\begin{tabular}{c|c|c|c|c|c}
\phantom{$\bullet$} &  &  &  &  &  \\ \hline
& \phantom{$\bullet$} &  &  &  &  \\ \hline
&  & $\bullet $ &  &  &  \\ \hline
&  &  & $\bullet $ & $\bullet $ &  \\ \hline
&  & $\bullet $ & $\bullet $ &  &  \\ \hline
&  &  &  &  & \phantom{$\bullet$}%
\end{tabular}%
\ \Rightarrow 
\begin{tabular}{c|c|c|c|c|c}
\phantom{$\bullet$} &  &  &  &  &  \\ \hline
& \phantom{$\bullet$} &  &  &  &  \\ \hline
&  &  & $\bullet $ &  &  \\ \hline
&  &  &  & $\bullet $ &  \\ \hline
&  & $\bullet $ & $\bullet $ & $\bullet $ &  \\ \hline
&  &  &  &  & \phantom{$\bullet$}%
\end{tabular}%
\end{eqnarray}%
It travels as a bishop in chess and in four moves it is shifted one cell
forward and down and its velocity is $v=1/4$. 

No object of the game of \textsc{Life} can move faster than a chess king --
one square per one move. Being called as \textit{the speed of light}, this
quantity is taken as the universal unit $c\equiv 1$. But the fact of '%
\textit{motion}' in discrete space-time of \textsc{Life} never obeys to
direct physical laws. The process of motion is based on \ the  intrinsic
nature of the cellular automata and the relevant velocity translation law
(if any) should be derived from it. The law which governs the relative
motion in the game of \textsc{Life} should be derived from the principle of
motion in discrete space-time.

The spaceships of the game of \textsc{Life} have direct resemblance with
chess pawns because they move at a discrete speed and strictly linear
(straight or diagonal) motion. The spaceship's volatility is originated from
the nature of cell unit which has only two possibilities -- either to
survive or to die. The process of motion follows from the death and birth
that periodically occurs in the cell society. 

The existence of spaceships is extremely important fact in the game of 
\textsc{Life} as it can be applied to construction of mathematical
(heuristic) models of the physical laws. Particularly, the observer flying
in the spaceship launches a 'bullet', and some \textit{external} party is
monitoring the events in his laboratory reference frame. The latter must
record the 'bullet' whose velocity is not coinciding with that measured in
the reference frame co-moving the 'astronaut'.

The recent attempt to derive special relativity from the cellular automata 
\cite{r3} inspired us to solve this vital problem: what is the formula of
velocity translation in the reference frame co-moving the spaceship? Either
it is the Galilean summation or the Lorentz relativistic formula? Although
there is no spaceship flying faster than $v=1/2$, it does not imply validity
of the plain Galilean summation. The attempts to derive special relativity
from the cellular automata \cite{r3} may hint to some answer. However, the
standard relativistic formula is also unmotivated because the motion in
discrete space-time of \textsc{Life} is not a continuous process. It is
composed of two primitive acts: elementary jump to the neighboring square
and staying at rest.

Is there any law to govern the relative motion between the moving reference
frames? Its importance cannot be overestimated. Its knowledge will turn the
cellular automata game of \textsc{Life} towards realities of the material
world.

\section{Motion on chessboard}

Consider a chess king on the chess board. Let the king always moves forward
and straight (never goes back or diagonal). Or, we can consider a pawn
instead (without the initial double step). Let after $P$ moves it appears
shifted in $n$ squares forward. We define its velocity as $v=n/P$. What has
happened within this time period? It is clear that the pawn was remaining at
rest during $P-n$ moves, while during the rest $n$ moves it was jumping one
square forward. In fact, the pawn has only two states: either to move or to
stay at rest. On the other hand, each move can be explained in terms of the
game of \textsc{Life} as follows: the old pawn dies and the new pawn appears
on a free square.

Now consider two pawns going on the same file of the chessboard to met each
other. 
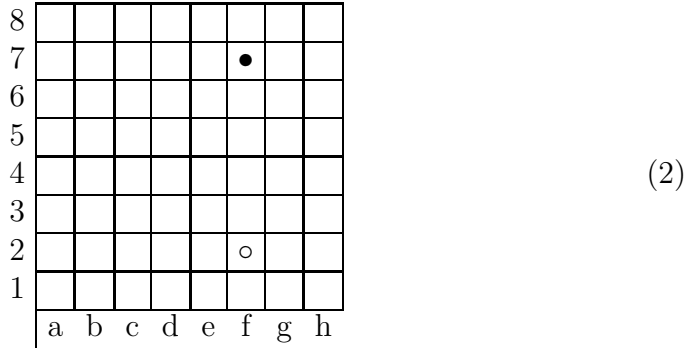
\begin{equation}
\begin{tabular}{l|c|c|c|c|c|c|c|c|}
\cline{2-9}
8 &  &  &  &  &  &  &  &  \\ \cline{2-9}\cline{9-9}
7 &  &  &  &  &  & $\bullet $ & \multicolumn{1}{|c|}{} &  \\ 
\cline{2-8}\cline{9-9}
6 &  &  &  &  &  &  & \multicolumn{1}{|c|}{} &  \\ \cline{2-8}\cline{9-9}
5 &  &  &  &  &  &  & \multicolumn{1}{|c|}{} &  \\ \cline{2-8}\cline{9-9}
4 &  &  &  &  &  &  & \multicolumn{1}{|c|}{} &  \\ \cline{2-8}\cline{9-9}
3 &  &  &  &  &  &  & \multicolumn{1}{|c|}{} &  \\ \cline{2-8}\cline{9-9}
2 &  &  &  &  &  & $\circ $ & \multicolumn{1}{|c|}{} &  \\ 
\cline{2-8}\cline{8-9}
1 &  &  &  &  &  & \multicolumn{1}{|c|}{} &  &  \\ \cline{2-9}
& \multicolumn{1}{c}{a} & \multicolumn{1}{c}{b} & \multicolumn{1}{c}{c} & 
\multicolumn{1}{c}{d} & \multicolumn{1}{c}{e} & \multicolumn{1}{c}{f} & 
\multicolumn{1}{c}{g} & \multicolumn{1}{c}{h}%
\end{tabular}
\label{chess}
\end{equation}
The white pawn \textbf{f2 }moves with the velocity $v_1$ and the black pawn%
\textbf{\ f7} moves with the velocity $v_2$ (it approaches just from north).
What is the net relative velocity? Checking it as a plain sum $v_1+v_2$ is
not correct because it may occur arbitrary high (the known Zeno's aporia) in
the series of consequent summation (of course, providing infinite
chessboard). But it is not so in the light of the chess rules because the
pawns do not move \textit{simultaneously}: the black pawn turns to move when
the white pawn stays at rest. If during $P$ moves the white pawn has passed $%
n$ squares, the black pawn has had a possibility to move exactly $P-n$
times. Therefore the net velocity cannot exceed $\left( n+P-n\right) /P$
that is the speed of light $c=1$.

\section{Relativistic translation of velocity}

Let an arbitrary spaceship moves at velocity $v_{1}$ and a hypothetical 
\textit{astronaut} is sitting inside. Let him have a gun and shot a \textit{%
bullet} flying as fast as $v_{2}$ (it is also a spaceship in the game of 
\textsc{Life} classification). The latter $v_{2}$ is the velocity with
respect to the reference frame co-moving the \textit{astronaut}. The
so-called 'guns' emitting spaceships are well known to the game of \textsc{%
Life} researches \cite{summers}.

Imagine that some \textit{external} party, an abstract  \textit{policeman}
is monitoring the events from his outpost based on the \textit{ground} or
the laboratory reference frame. First, he records the \textit{astronaut}
flying in his spaceship with the velocity $v_{1}$. Second, he records the 
\textit{bullet} whose velocity is $v_{12}$ and we shall calculate it now

Again the naive summation $v_{1}+v_{2}$ may exceed the speed of light. The
game of \textsc{Life} spaceship is a complex object composed of many cells
but as a single pawn, considered above, it is involved in the motion
consisting in $n$ pure jumps and $P-n$ stays because after $P$ moves the
spaceship is replaced into $n$ squares. In the reference frame of the 
\textit{astronaut }in the spaceship his proper time is staying at rest and
it corresponds to $P-n$ moves. During his proper time he records the \textit{%
bullet} passing the distance $s_{2}$ at its native velocity\thinspace\ $%
v_{2}=s_{2}/\left( P-n\right) $.

The \textit{policeman} records the spaceship passed $n$ cells (at the
velocity $v_{1}=n/P$) and the total path of the \textit{bullet} is recorded
as\thinspace\ $s_{1}+s_{2}=v_{1}P+v_{2}\left( P-n\right) $. Therefore, the 
\textit{policemen} has finally recorded the \textit{bullet} flying with the
velocity 
\begin{equation}
v_{12}=v_{1}+v_{2}-v_{1}v_{2}  \label{v}
\end{equation}%
It is a translation formula acting on the discrete cellular board of the
game \textsc{Life}.

Since $v_{1}<1$ and $v_{2}<1$, the total velocity cannot exceed the speed of
light. Particularly, light always propagates at its genuine speed $v_{2}=c=1$
because $v_{12}=1$ without regard of the velocity of the reference frame $%
v_{1}$. Another sample: at $v_{1}=v_{2}=1/2$ the relative velocity is $%
v_{12}=3/4$ (but not the plain sum!).

Now let us have some recreation and solve the real situation. Imagine, you
are rich enough to purchase any available vehicle and you have chosen a
spaceship moving at velocity $2/5$. Of course, you have a gun \cite{summers}
which can emit e.g. \textit{lightweight} spaceships (SWSS) with velocity $%
v_{2}=1/2$ \cite{koenig}. However... your are arrested because your gun is
expected as dangerous: the bullets faster than $3/4$ of the speed of light
are forbidden. Meanwhile, when flying in your spaceship you can shoot
bullets at $v_{12}=v_{1}+v_{2}=9/10$. The relativistic composition of
velocities yields $v_{12}=(v_{1}+v_{2})/(1+v_{1}v_{2})=3/4$. At any rate you
cannot avoid jail. What shall you do? Of course, you smile and refer to the
formula (\ref{v}) resulting to $v_{12}=7/10$. This cellular world is
governed by the mathematical rules and laws and it is not derived from the
real world... 

Let us put our \textit{astronaut} in a luxury puffer-class spaceship (or
briefly \textit{puffer} \cite{summers} according to the researches of the
game of \textsc{Life}). It is not a spaceship in the strict sense but also a
flying vehicle which, however, emits gliders. Glider (\ref{g}).belongs to
the category of spaceships moving oblique course but it is not a tedious
problem to calculate the composition of velocities for this \textit{bullet}.
Again, consider a manned puffer flying at the velocity $v_{1}$ along the
axis $x$. The proper time of the \textit{astronaut} is $P\left(
1-v_{1}\right) $. Let a \textit{bullet} is shot from the spaceship at some
angle $\psi $ with respect to the axis $x$ and this angle is defined as 
\begin{equation}
\tan \psi =v_{2}^{y}/v_{2}^{x}
\end{equation}%
where $v_{2}^{x}$ and $v_{2}^{y}$ is the velocity of the \textit{bullet}
along the axes $x$ and $y$ respectively. Of course, it is possible to
introduce some nominal velocity $v_{2}$ so that 
\begin{equation}
v_{2}^{x}=v_{2}\cos \psi \qquad v_{2}^{y}=v_{2}\sin \psi 
\end{equation}%
however, this quantity has not real physical meaning of motion at the angle $%
\psi $ with the velocity $v_{2}$. The Pythagoras theorem is not applied
here. Indeed, a glider (\ref{g}) is walking along the diagonal. Its path is
one square forward and one square down in four moves, hence, yielding
velocity \ $v_{x}=v_{y}=1/4$ and its velocity along the diagonal is also $1/4
$ \ rather than $\sqrt{2}/4$.

The \textit{astronaut} will see the \textit{bullet} flying the distance $%
s_{2}^{x}=P\left( 1-v_{1}\right) \,v_{2}^{x}$ along the axis $x$ and the
distance $s_{2}^{y}=P\left( 1-v_{1}\right) \,v_{2}^{y}$ \ along the axis $y$
(which is orthogonal to the axis $x$).

The relevant path of the \textit{bullet} in the \textit{earth} reference
frame along the axis $x$ is extended by the path $s_{1}=n$ of the carrier
spaceship. The summary boost of the \textit{bullet} along the axis $x$ will
be $s^{x}=n+P\left( 1-v_{1}\right) \,v_{2}\cos \psi $ \ implying that its
velocity is 
\begin{equation}
v_{12}^{x}=v_{1}+\left( 1-v_{1}\right) \,v_{2}^{x}  \label{vv}
\end{equation}%
meanwhile, the summary replacement along the axis $y$ has not become longer
than $s^{y}=s_{2}^{y}$ and we write 

\begin{equation}
v_{12}^{y}=\left( 1-v_{1}\right) v_{2}^{y}  \label{yy}
\end{equation}%
Note that this component does not coincide with $v_{2}^{y}$ .

The 'bullet' flies in the direction defined in the \textit{earth} reference
frame by the angle $\chi $ so 
\begin{equation}
\tan \chi =\frac{\left( 1-v_{1}\right) \,v_{2}\sin \psi }{v_{1}+\left(
1-v_{1}\right) \,v_{2}\cos \psi }  \label{angle}
\end{equation}

Again, the total velocity cannot exceed the speed of light $c=1$ and the
magnitude of (\ref{vv})-(\ref{yy}) is always smaller than the value (\ref{v}%
) to which it is reduced when the \textit{bullet} is shot parallel to the
course of 'carrier' spaceship ($\psi =0$) and the deviation from the special
relativity 
\begin{equation}
\Delta =v_{1}v_{2}\frac{1-v_{1}-v_{2}+v_{1}v_{2}}{1+v_{1}v_{2}}
\end{equation}%
never exceeds $0.05$.

\section{Sample}

There exist \textit{puffers} \cite{summers} consisting of about 1000 cells
and drifting at a velocity $v_{1}=1/4$. Let our \textit{astronaut} drive one
of them.  There are many orthogonal spaceships (i.e. $\psi =0$) with
velocity $v_{2}=1/3$. Some of them are small and good enough to be a \textit{%
bullet }which the \textit{astronaut can }shot at the velocity

\begin{equation}
v_{2}^{x}=0\qquad v_{2}^{y}=1/3  \label{em0}
\end{equation}%
in the co-moving reference frame. 

What shall see the \textit{earth }observer in the \textit{earth} reference
frame? According to formulas (\ref{vv})-(\ref{angle}) the external observer
will record the object whose velocity is%
\begin{equation}
v_{12}^{x}=-1/4\qquad v_{12}^{y}=1/4  \label{em1}
\end{equation}%
whose course is inclined at the angle 
\begin{equation}
\chi =135^{\circ }  \label{em2}
\end{equation}%
with respect to the course of the carrier spaceship. In fact, among the 
\textit{puffers} \cite{summers} drifting at a velocity $v_{1}=1/4$ there are
several patterns which emit gliders exactly in the prescribed direction (\ref%
{em1})-(\ref{em2}) . 

What about the \textit{bullet} spaceship (\ref{em0})? We do not observe it
explicitly if we sit on the \textit{earth}. How we can grasp it? How we can
interpret it? It can be recognized as a definite texture within the \textit{%
puffer} pattern, an \textit{embryo} carried by the carrier spaceship until
it is pushed out as a glider of the \textit{earth} reference frame. The
virtual \textit{embryo} is moving at the velocity (\ref{em0}) within the
bosom of the carrier \textit{puffer} and its embryonal motion is resulted
from the interaction of cells. As soon as the \textit{embryo} is born, it
leaves the womb of its 'mother' spaceship and goes in the outer space where
it becomes glider visible explicitly even to the external observers.\bigskip 

\section{Conclusion}

The velocity translation law of the game of \textsc{Life} is discovered (\ref%
{vv})-(\ref{angle}). It differs from the relevant Lorentz formula of special
relativity. Indeed, the motion in the game of \textsc{Life} is resulted from
the given rules of the cellular automata and has no explicit physical
nature. Besides there is evident anisotropy of space in vertical and
horizontal direction, while the orthogonal and diagonal motion occurs at the
same velocity (there is no Pythagoras theorem).

The formulas (\ref{vv})-(\ref{angle}) allow to operate with the moving
reference frames. Since before this possibility was not expected in the game
of \textsc{Life}. Now we have the approach to consider many interesting
ideas, e.g. 'the light communication' or the telegraph based on the
spaceships (the pilot emits signals propagating at various speed etc). The
quantitative expression of the velocity translation law (\ref{vv})-(\ref%
{angle}) opens new horizon for formulation of the 'relativistic' problems in
the game of \textsc{Life}.

\end{document}